# Period-doubling cascades of a Silnikov equation

Keying Guan† and Beiye Feng‡

†Science College, Beijing Jiaotong University, Email: keying.guan@gmail.com

‡Institute of Applied Mathematics, Academia Sinica, Email: fby@amss.ac.cn

Abstract: Based on numerical results of a Silnikov equation, three period-doubling cascades, corresponding respectively to three different characters of the rotation number of a limit closed orbit, are studied, and the Feigenbaum constant is used successfully in the estimation of the critical parameter values for the period-doubling bifurcation. The conceptions of separaror and pseudo-attractor are also introduced in the discussion part of this paper.

1. Introduction

In the paper [1], the author considered the dynamical system (a particular Silnikov equation)

$$\begin{cases} \frac{dx}{dt} = P(x,y,z) = y \\ \frac{dy}{dt} = Q(x,y,z) = z \\ \frac{dz}{dx} = R(x,y,z) = x^3 - a^2 x - y - b z \end{cases} \quad (1)$$

and found that this system has two kinds of non-trivial attractors, one kind is with complex structures（maybe fractal structure）, and the other kind is some spatial limit closed orbits（including the limit cycles）. Some bifurcation phenomena about the limit closed orbits had also been discussed at that paper.

The authors have noticed that the results mentioned above are in a good consistency with the theory of the period-doubling cascades which is studied systematically for the abstract dynamical systems by Evelyn Sander and James A. Yorke (ref. [2],[3]).

Though, the abstract dynamical system may largely represent the fundamental properties of the autonomous ordinary differential equation system, but each concrete given ODE system has still its own particular bifurcation phenomena. It is proven in this paper that the character of the rotation number is an important new conception for the determination of the series of period-doubling cascades. This new conception comes just from our research to the particular ODE system (1). Besides, it is usually very difficult to find a concrete ODE system to realize all of the phenomena predicted by the theory of the period-doubling cascades of the abstract dynamical system.

Therefore, the study of the period-doubling cascades for a given ODE system has still its own particular significance.

## 2. Some fundamental properties and facts of the system (1)

The ODE system (1) has a series of nice properties for the study of the period-doubling cascades:

$P_1$. The system (1) is symmetrical about the origin. It has three equilibrium points, $p_0 = (0,0,0)$, $p_1 = (-a, 0,0)$ and $p_2 = (a, 0,0)>$. When $a = 1$, $0 \leq b < 1$, all of the three equilibrium points are saddle-focus type. Concretely, the equilibrium point $p_0$ has a one-dimensional stable manifold $M_{1d}^s(p_0)$ and a two-dimensional unstable manifold $M_{2d}^u(p_0)$ that is corresponding to a pair of conjugate eigenvalues with positive real part, while both equilibrium points $p_1$ and $p_2$ have, respectively, a one-dimensional unstable manifold $M_{1d}^u(p_i)$ and a two-dimensional stable manifold $M_{2d}^s(p_i)$ (i=1,2) that corresponds to a pair of conjugate eigenvalues with negative real part.

$P_2$. When $a = 1$, $b_0 \leq b < 1$ (the recent numerical results show $b_0 \approx 0.3121$), then the closed set

$$A_{class} = \overline{M_{2d}^u(p_0) \backslash (M_{2d}^u(p_0) \cup \{p_0, p_1, p_2\})} \quad (2)$$

is bounded. The phase space of system (1) has an invariant sub-manifolod $\mathfrak{M}$, which is the region between the two-dimensional stable manifolds $M_{2d}^s(p_1)$ and $M_{2d}^s(p_2)$, and all of the integral orbits in $\mathfrak{M}$ approach $A_{class}$.

The following figures (in the case, $a = 1$, $b = 0.315$) are used to show these conceptions and facts intuitively.

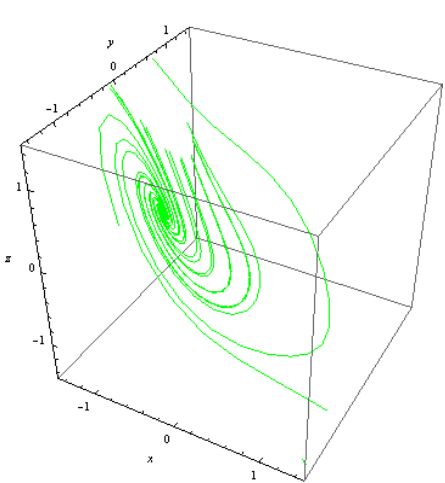 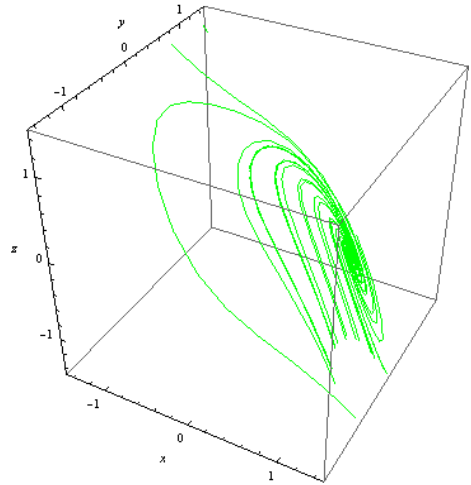

Figure 1. $M_{2d}^s(p_1)$    Figure 2. $M_{2d}^s(p_2)$

Figures 1 and 2 shows the two-dimensional stable manifolds $M^s_{2d}(p_1)$ and $M^s_{2d}(p_2)$ respectively, where each manifold is expanded with its integral curves (colored green).

Figure 3 shows the set $A_{class}$ (colored purple), which is a nontrivial attractor with complex structure. Figure 4 shows an integral curve (colored blue) approaches the attractor $A_{class}$.

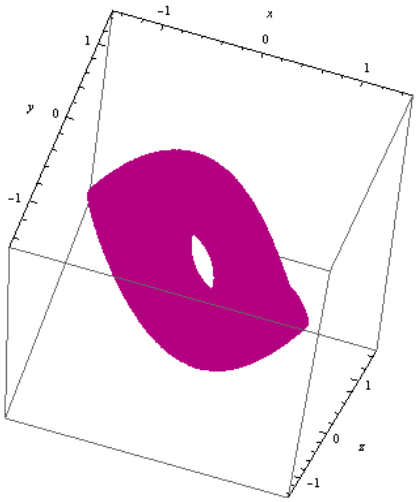 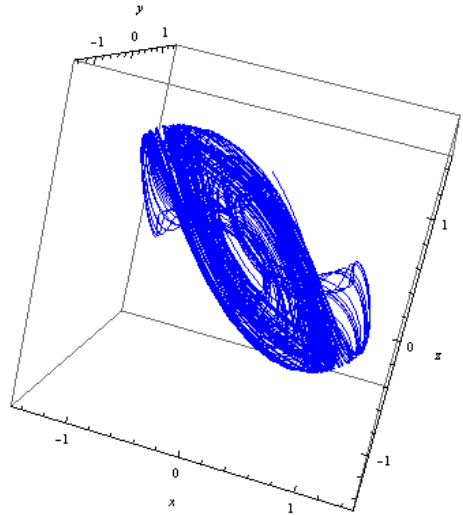

Figure 3.  $A_{class}$          Figure 4. An integral curve approaching $A_{class}$

Figure 5 and 6 shows the geometrical relations between the above-mention objects from two different viewing angles.

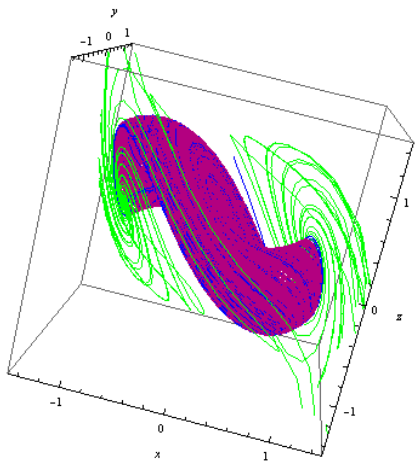 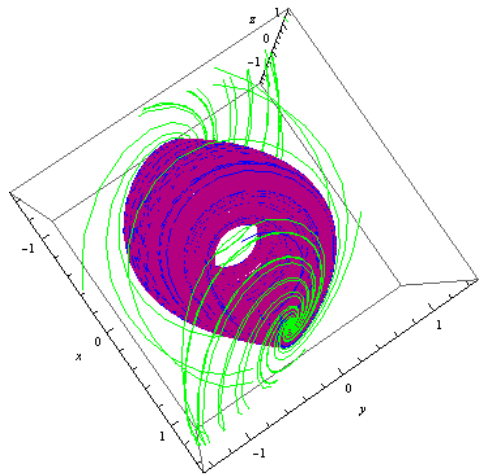

Figure 5. From one viewing angle          Figure 6. From another viewing angle

It is easy to see this invariant manifold $\mathfrak{M}$ is quite similar to the one used in the article [3].

Besides the general properties $P_1$ and $P_2$, the attractor may be a spatial closed orbit with a certain rotation number about the corresponding branch of the one-dimensional stable manifold $M_{1d}^s(p_0)$ (ref. [1]). The rotation numbers of the system (1), found numerically by now, are in three different series:

*Series 1.*  $1, 2, 4, \cdots, 2^n, \cdots,$
*Series 2.*  $13, 2 \times 13, 4 \times 13, \cdots, 2^n \times 13, \cdots,$
*Series 3.*  $3, 2 \times 3, 4 \times 3, \cdots, 2^n \times 3, \cdots$

It is easy to see that the three numbers, $1, 13$ and $3$ are the character of the three series of rotation numbers respectively. Clearly, it is easy to give a general definition to the topological character of rotation number for a given limit closed orbit in a strict way. However, in order to concentrate to the study of the relevant phenomena first, this paper would like to neglect the strict and general definition.

When the corresponding parameter $b$ varies, there would be a series of critical values (bifurcation points) for the sudden change (doubling) of the rotation number. If the rotation number of a closed orbit is doubling, the periodic of the closed orbit is also doubling obviously. So the above-mentioned three series must correspond to three cascades of period-doubling.

The numerical results show that the three series of critical values are distributed in three different intervals of the parameter $b$. This fact shows the importance of the topological character of rotation number of the closed limit orbit.

The next section will introduce the concrete numerical results for these cascades.

### 3. Period-doubling doubling and Feigenbaum constant

3.1 Period-doubling bifurcation for *series 1*.

If $b_{1,inf} < b < 1$, the bifurcation points are distributed as

$$b_{1,inf} < \cdots < b_{1,n} < \cdots < b_{1,3} < b_{1,2} < b_{1,1} < b_{1,0} < 1 \qquad (3)$$

such that,

(1.1) when $b_{1,0} \leq b < 1$, the system (1) has a unique limit cycle of rotation number 1;
(1.2) when $b_{1,n} \leq b < b_{1,n-1}$, the system (1) has a pair of spatial limit closed orbits of rotation number $2^n$, they are symmetric about the origin;
(1.3) $b_{1,inf} = \lim_{n \to \infty} b_{1,n}$.

The numerical results show that, $b_{1,0} \approx 0.4892$, $b_{1,1} \approx 0.3992$, $b_{1,2} \approx 0.3829$, $b_{1,3} \approx 0.3794$, and $b_{1,inf} \approx 0.3750$.

Figure 7 shows the unique limit cycle of (1) in the case $b = 0.4892$.

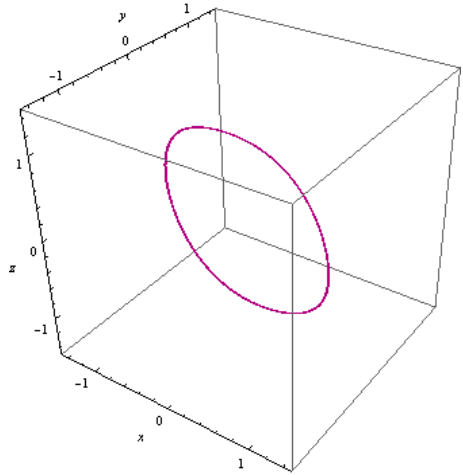

Figure 7. Single limit cycle

Figure 8 shows the pair of symmetric limit cycles in the case $b_{1,0} = 0.4891$, which are separated just from a single limit cycle.

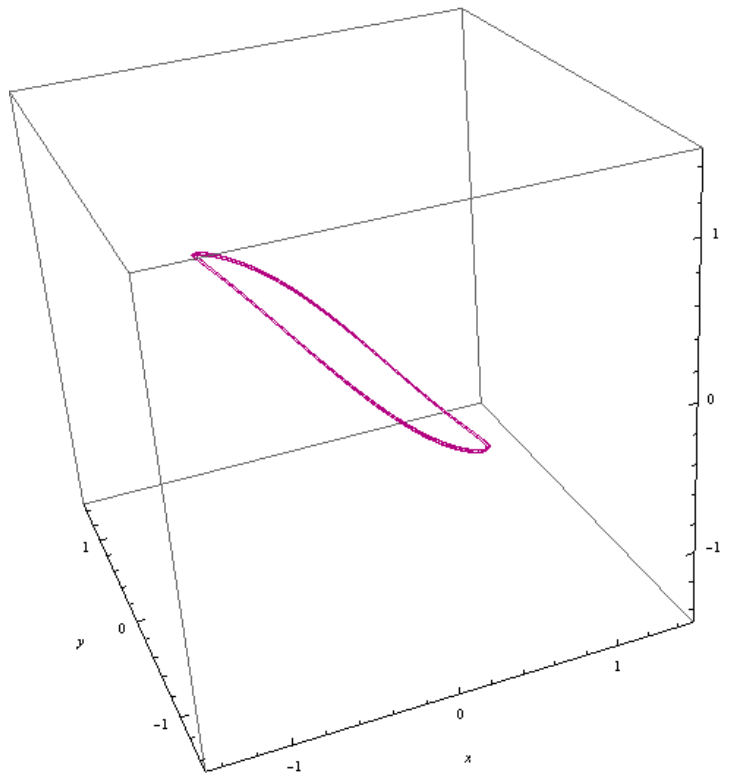

Figure 8. A pair of symmetric limit cycles

Figure 9 shows the pair of symmetrical limit cycles in the case $b = 0.3993$. Figure 10 shows the pair of symmetrical limit orbits with rotation number 2 after the period-doubling in the case $b = 0.3991$.

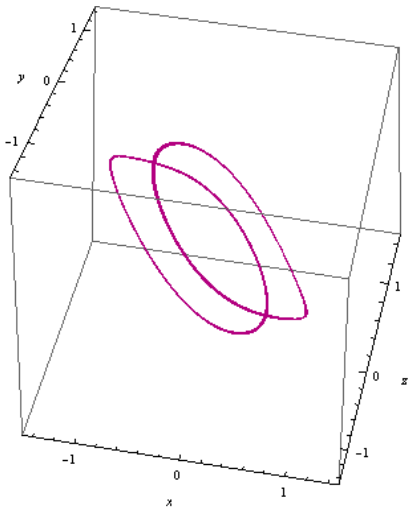
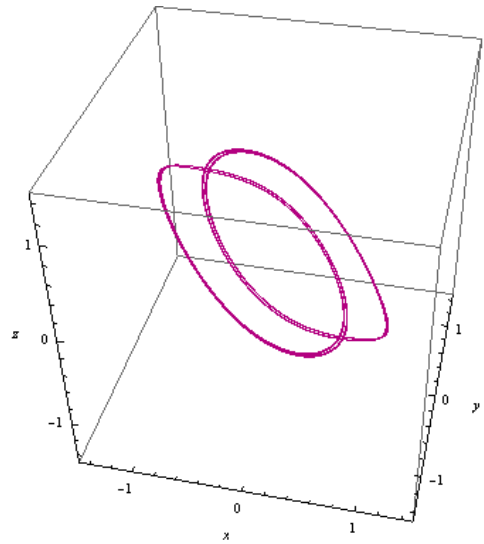

Figure 9. The pair of limit cycles before period-doubling

Figure 10. The pair of limit closed orbits after period-doubling

Figure 11 shows one of the pair of limit closed orbits of rotation number 2 in the case $b = 0.3830$, just before period-doubling. Figure 12 shows the limit closed orbit of the rotation number 4 in the case $b = 0.3828$, after the period-doubling.

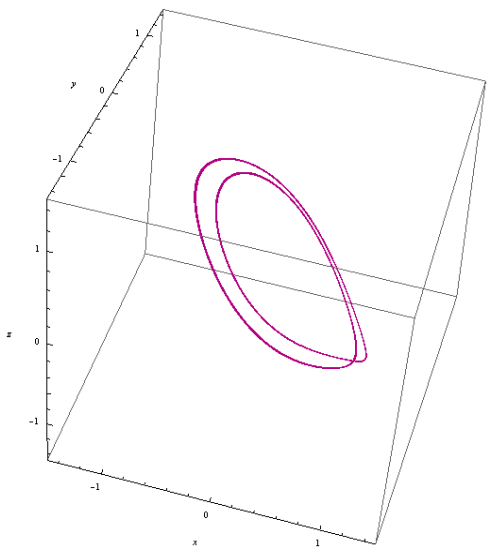
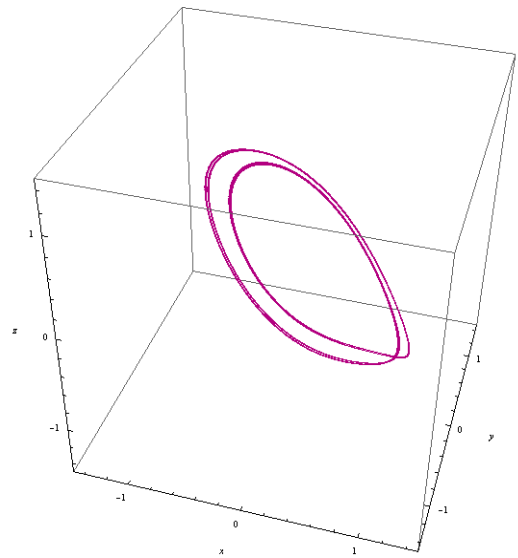

Figure 11. One limit closed orbit of rotation number 2 before the period-doubling

Figure 12. The limit closed orbit of rotation number 4 after the period-doubling

Figure 13 shows one of the pair of limit closed orbits with rotation number 4 in the case $b = 0.3795$, just before period-doubling. Figure 14 shows one limit closed orbit with the rotation number 8 in the case $b = 0.3793$, after the period-doubling. Figure 15 shows an enlarged part of the Figure 14, so that the rotation number 8 can be seen clearer.

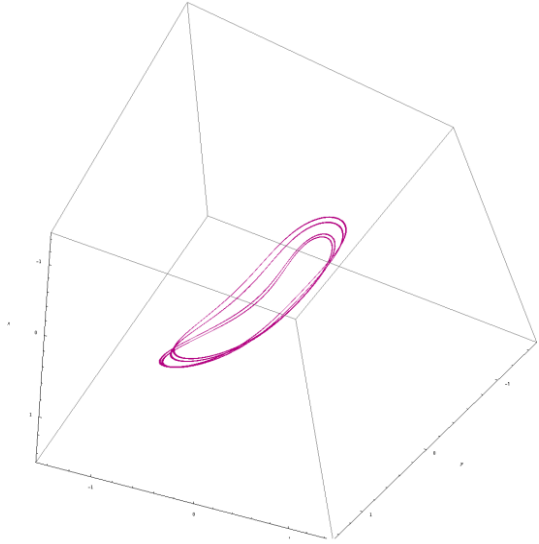

Figure 13. One limit closed orbit of rotation number 4

before the period-doubling

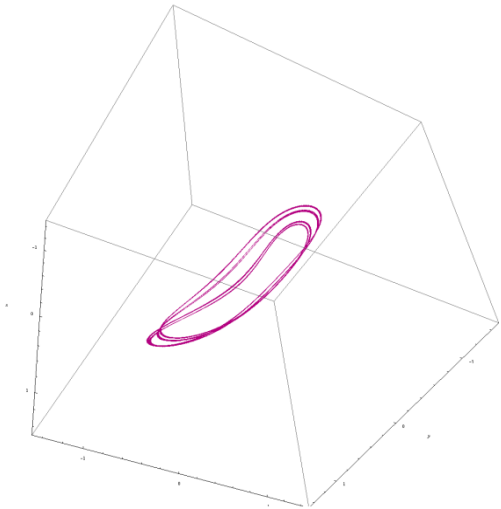

Figure 14. The limit closed orbit of rotation number 8

after the period-doubling

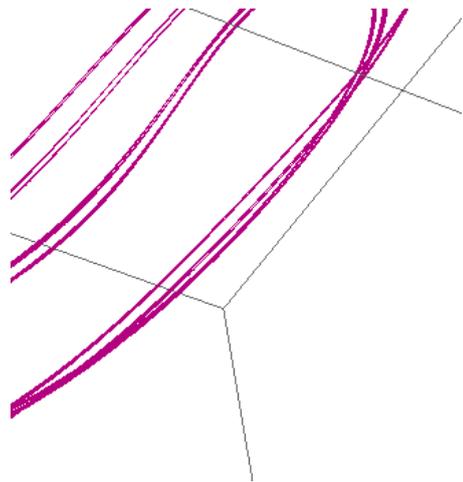

Figure 15. Enlarged part of Figure 14

Obvious, it is a difficult issue to estimate exactly the rotation number for the numerical method when the number is large. The estimation of $b_{1,inf}$ is more difficult.

However, the dynamical system theory may provide some help for the estimation of some critical values and $b_{1,inf}$. It is well known that, in the dynamical system theory, there is a famous Feigenbaum constant, $\delta = 4.669201\cdots$. The following is a brief introduction in Wikipedia:

*The first Feigenbaum constant is the limiting ratio of each bifurcation interval to the next between every period doubling, of a one-parameter map*

$$x_{i+1} = f(x_i)$$

where $f(x)$ is a function parameterized by the bifurcation parameter a.

It is given by the limit:

$$\delta = \lim_{n \to \infty} \frac{a_{n-1} - a_{n-2}}{a_n - a_{n-1}} = 4.669201609 \cdots$$

where $a_n$ are discrete values of a a at the $n$th period doubling.

Though this constant is obtained in the case that the mapping $f(x)$ is for the one dimensional veriable $x$, but the constant is believed to be universal.  So, it is a good chance here to test the universality of this constant in the following way:

Let the bifurcation value $b_{1,n}$ be treated as $a_n$, and assume that the equality

$$\frac{a_{n-1} - a_{n-2}}{a_n - a_{n-1}} = \delta \tag{4}$$

is satisfied approximately for $b_{1,2}$, $b_{1,3}, \cdots, b_{1,n}, \cdots$. Since we have numerically obtained $b_{1,2} \approx 0.3829$ and $b_{1,3} \approx 0.3793$, then from (4), it is easy to get the following estimations

$$b_{1,4} \approx b_{1,3} - \frac{b_{1,2} - b_{1,3}}{\delta} \approx 0.3785 \tag{5}$$

and

$$b_{1,inf} \approx b_{1,2} - \frac{b_{1,2} - b_{1,3}}{(1 - \frac{1}{\delta})} \approx 0.3783 \tag{6}$$

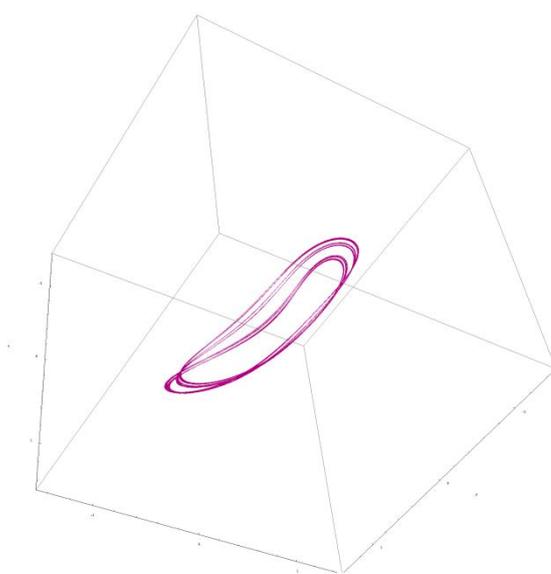
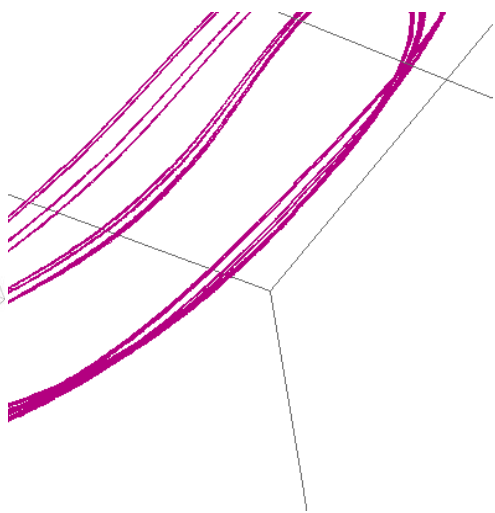

Figure 16. The limit closed orbit of rotation number 8 in the case $b = 0.3788$

Figure 17. Enlarged part of Figure 16

Based on the direct numerical calculations, Figure 16 and 17 shows that the rotation number of the limit closed orbit is 8 in the case $b = 0.3788$

Figure 18 and 19 shows clearly that the that the rotation number of the limit closed orbit is 16 after a period-doubling bifurcation in the case $b = 0.3786$. So, the direct numerical calculation of the system (1) gives $b_{1,4} \approx 0.3787$. This result shows that (5) gives a very good approximation for the bifurcation value $b_{1,4}$.

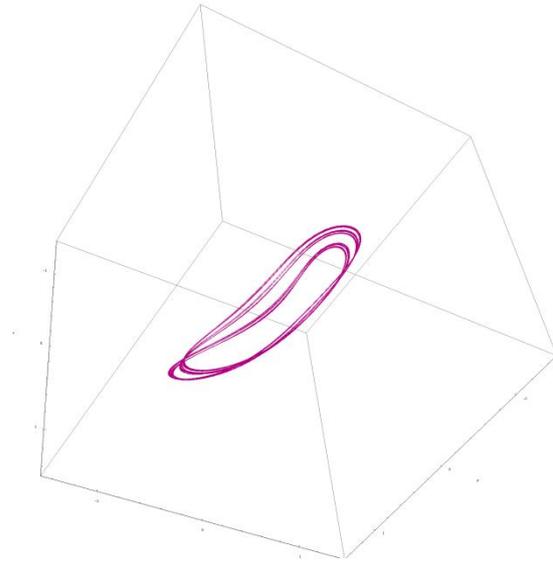
Figure 18. The limit closed orbit of rotation number 16 after a period-doubling in the case $b = 0.3786$

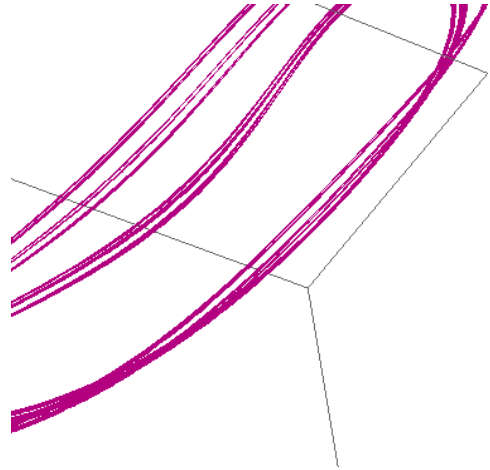
Figure 19. Enlarged part of Figure 18

Figure 20 and 21 shows that the set $A_{class}$ is formed with limit closed orbits in the case $b = 0.3785$

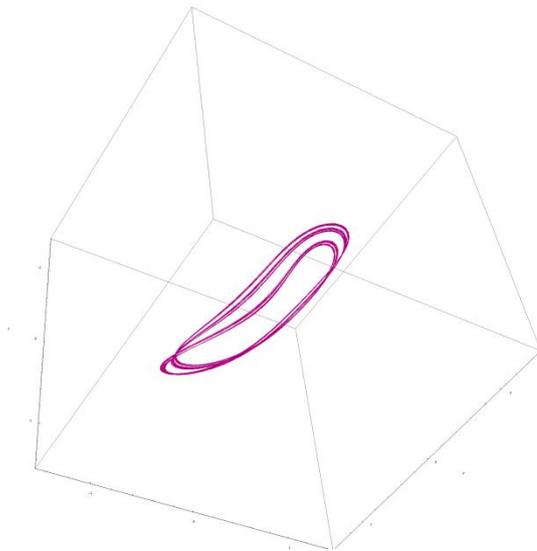
Figure 20. The set $A_{class}$ in the case $b = 0.3783$

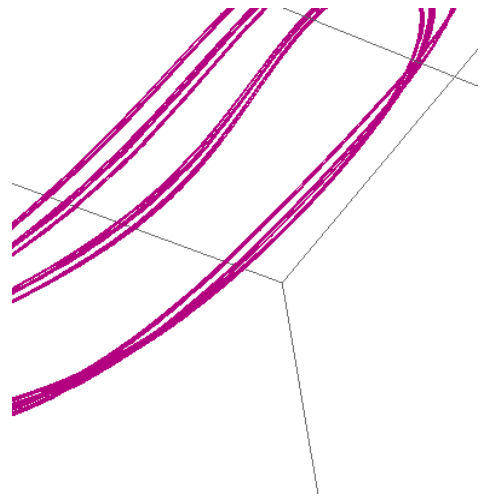
Figure 21. Enlarged part of Figure

Figure 22 and 23 suggests that the set $A_{class}$ should not be formed with spatial limit orbits in the case $b = 0.3783$.

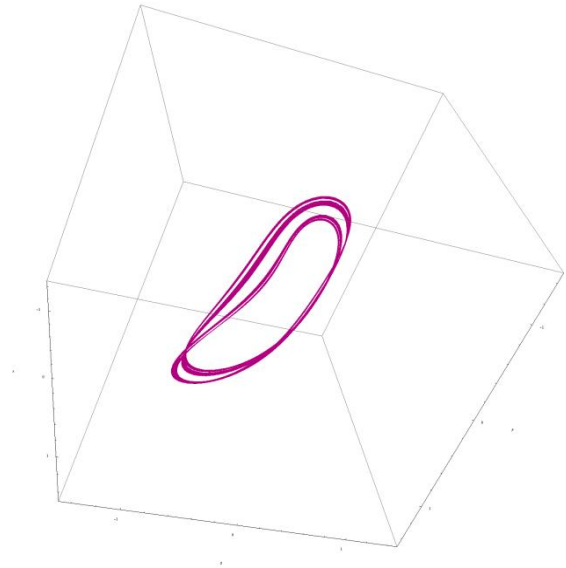 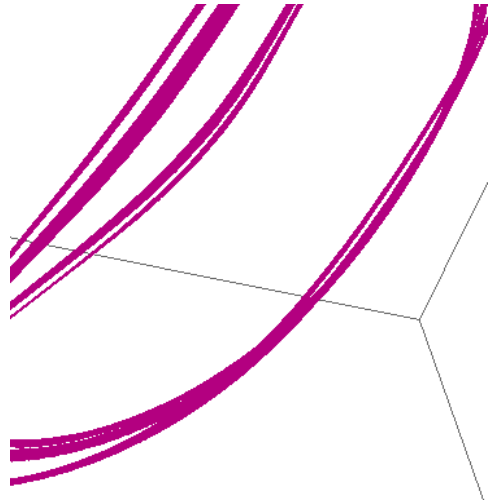

Figure 22. The set $A_{class}$ in the case $b = 0.3783$     Figure 23. Enlarged part of Figure 22

This fact shows that (6) gives a surprisingly good estimation for $b_{1,inf}$.

3.2 Period-doubling bifurcation for *series 2.*

When the parameter $b$ is in the closed interval $[b_{13}, b_{1,inf}]$, the set $A_{class}$ has a complex structure. This paper would not try to explore the complex structure of the set $A_{class}$ in this parameter interval. Numerical test shows that $b_{13} \approx 0.33415$.

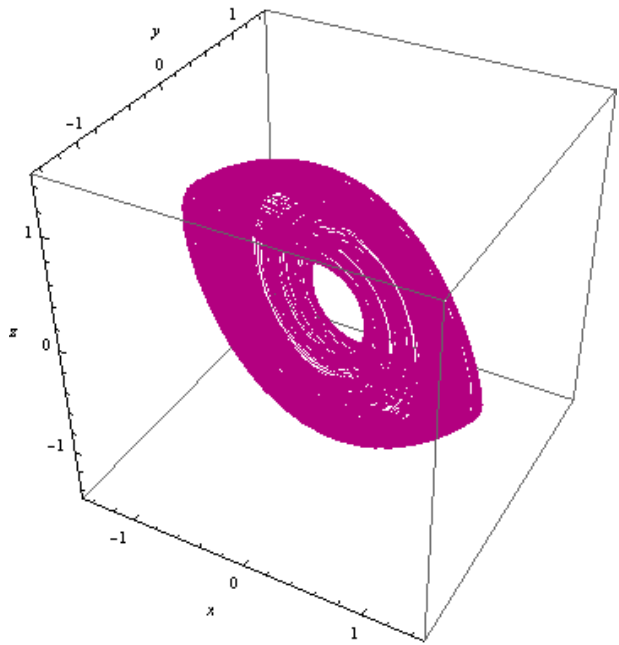

Figure 24. The set $A_{class}$ in the case $b = 0.3342$

In fact, Figure 24 shows that the set $A_{class}$ has still complex structure in the case $b = 0.3342$, and Figure 25 shows the set $A_{class}$ is a single spatial limit closed orbit with rotation number 13 in the case $b = 0.3341$. This orbit is symmetrical about the origin.

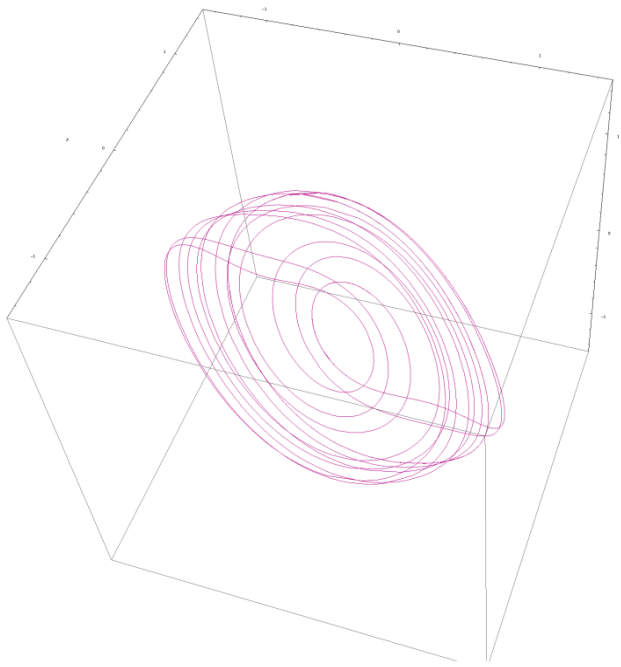

Figure 25. The limit closed orbit with rotation number 13

in the case $b = 0.3341$

Assume that, when the parameter $b$ is in the open interval, $(b_{13,inf}, b_{13})$, the set $A_{class}$ is formed with a single limit closed orbit of rotation number $2^n \times 13$ ($n = 0,1,2,\cdots$), and the limit closed orbit is symmetric about the origin. There should be a series critical values, $b_{13,n}$, $n = 0,1,2,\cdots$, distributed as follows

$$b_{13,inf} < \cdots < b_{13,n} < \cdots < b_{13,1} < b_{13,0} < b_{13}$$

such that

(2.1)  when $b_{13,0} \leq b < b_{13}$, the system (1) has a unique limit closed orbit of rotation number 13;

(2.2)  when $b_{13,n} \leq b < b_{13,n-1}$, the system (1) has a pair of spatial limit closed orbits of rotation number $2^n \times 13$;

(2.3)  $b_{1,inf} = \lim_{n \to \infty} b_{13,n}$.

Figure 26 shows that the rotation number of the limit closed orbit is 13 in the case $b = 0.3332$

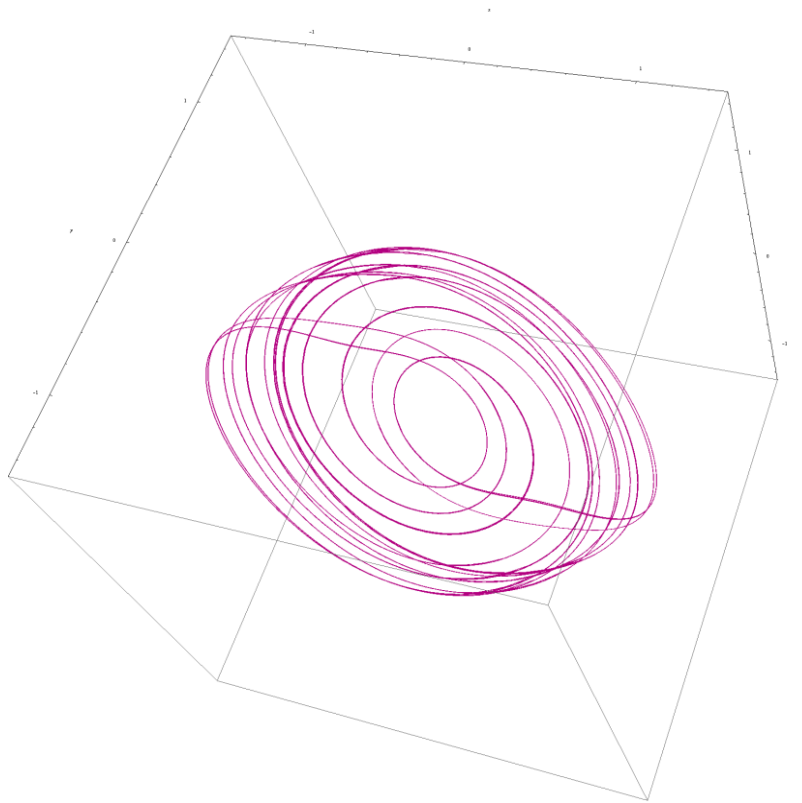

Figure 26. The limit closed orbit with rotation number 13

in the case $b = 0.33332$

Figure 27, especially Figure 28 (the enlarged part of Figure 27), shows the rotation number of the limit closed orbit has been doubled to 26 in the case $b = 0.33330$.

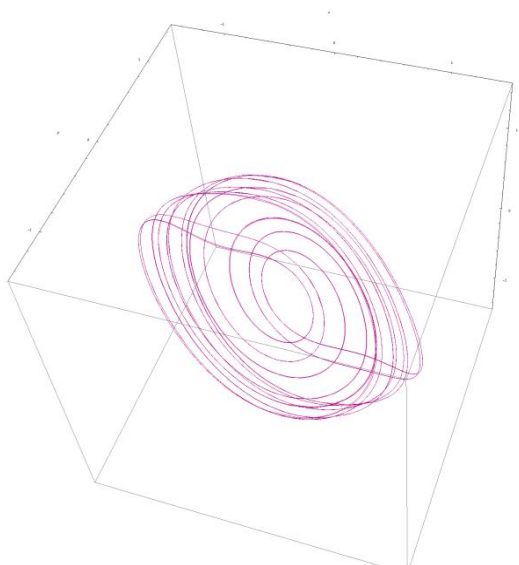 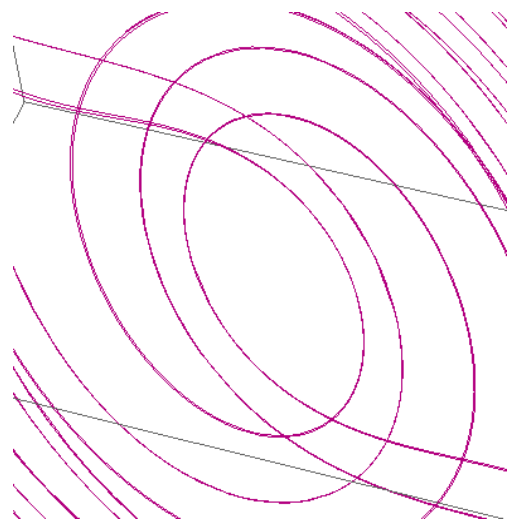

Figure 27.　The limit closed orbits　　　　　Figure 28. Enlarged part of Figure 27

in the case $b = 0.33330$

The above facts show that $b_{13,0} \approx 0.33331$

Figure 29 and 30 shows the rotation number of the limit closed orbit has been doubled to 56 in the case $b = 0.3333298$.

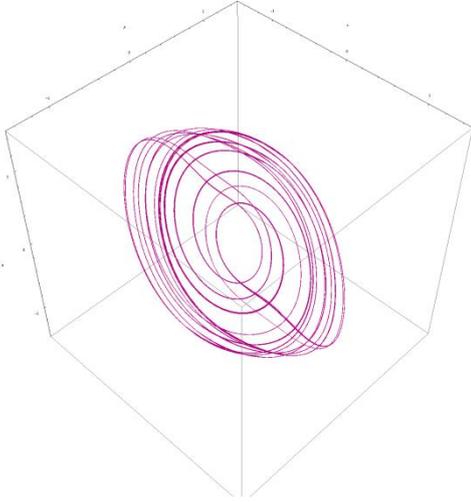 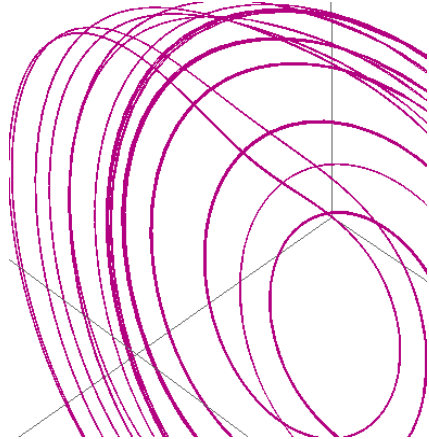

Figure 29.　The limit closed orbits　　　　　　Figure 30. Enlarged part of Figure 29
in the case $b = 0.333298$

From Figure 28 and 30, it is easy to see that $b_{13,1} \approx 0.333299$.

Since $b_{13,0} - b_{13,1} < 0.00005$，the further numerical calculations and checks for $b_{13,2}$, $b_{13,3}$, ⋯, become difficult more and more.  Instead, the value of $b_{13,inf}$ may be obtained by the estimation

$$b_{13,inf} \approx b_{13,0} - \frac{b_{13,0} - b_{13,1}}{(1 - \frac{1}{\delta})} \approx 0.333295 \qquad (7)$$

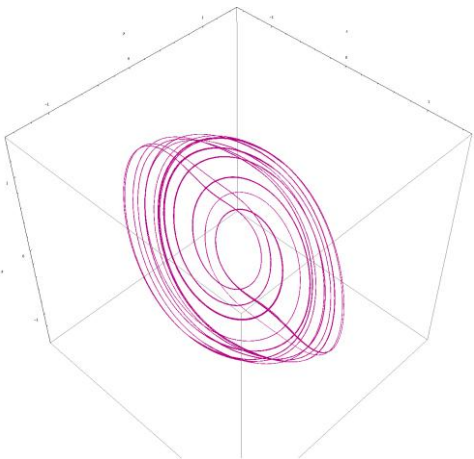 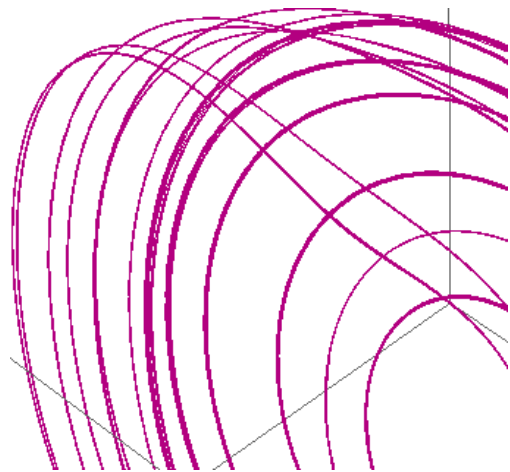

Figure 31.　The limit closed orbits　　　　　　Figure 32. Enlarged part of Figure 31
in the case $b = 0.333296$

Based on the numerical calculation, Figure 31 and 32 shows that, the set $A_{class}$ is still formed with a limit closed orbits in the case $b = 0.333296$, though its rotation number is difficult to be counted.

However, the numerical calculation suggests that the set $A_{class}$ should not be formed with a limit closed orbits.  See Figure 33 and 34.

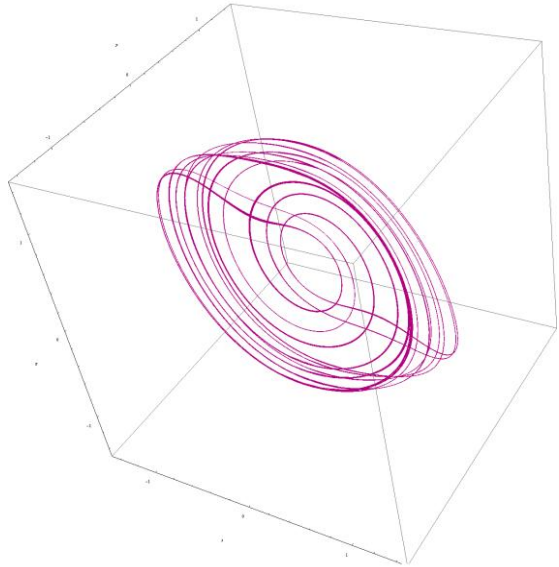
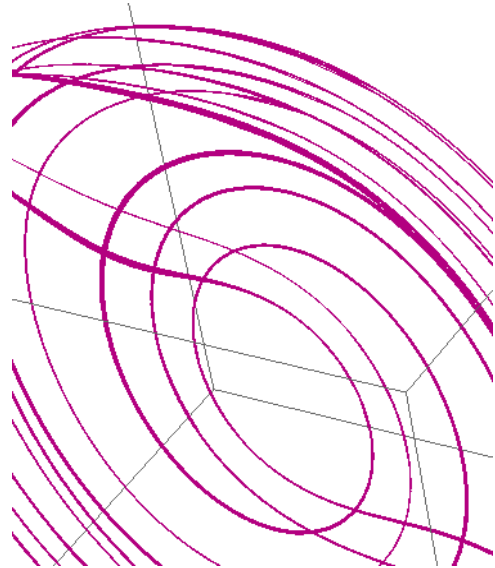

Figure 33.  A branch of $A_{class}$ in the case $b = 0.333295$

Figure 34. Enlarged part of Figure 33

So, the Feigenbaum constant is proven to be very effective for the estimation of $b_{13,inf}$.

3.3  Period-doubling bifurcation of *series 3.*

The numerical results show that there is an open interval $(b_{3,inf}, b_3)$ for the parameter $b$, such that, the system (1) has a pair of spatial limit closed orbits with rotation number, $2^n \times 3$, $n = 0, 1, 2, \cdots$.

The numerical test suggests that $b_3 \approx 0.32375$, since the set $A_{class}$ is with complex structure when $b = 0.3238$ (see Figure 35), and the set $A_{class}$ is formed with a pair of limit closed orbits with rotation number 3 when $b = 0.3237$ (see Figure 36，37).

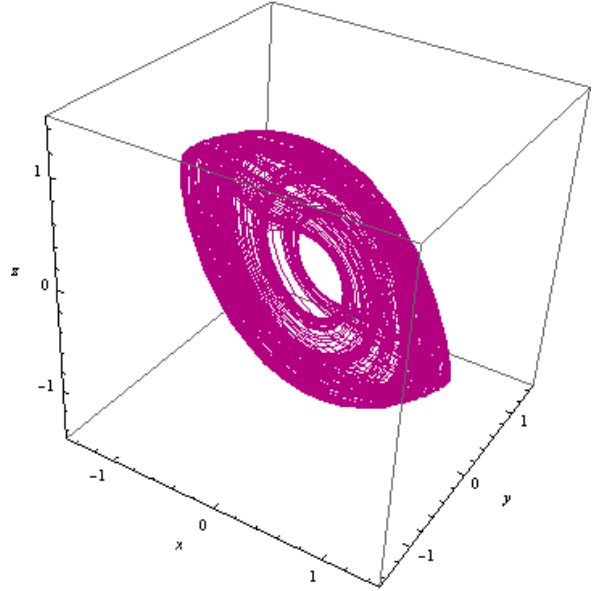

Figure 35. The set $A_{class}$ in the case $b = 0.3238$

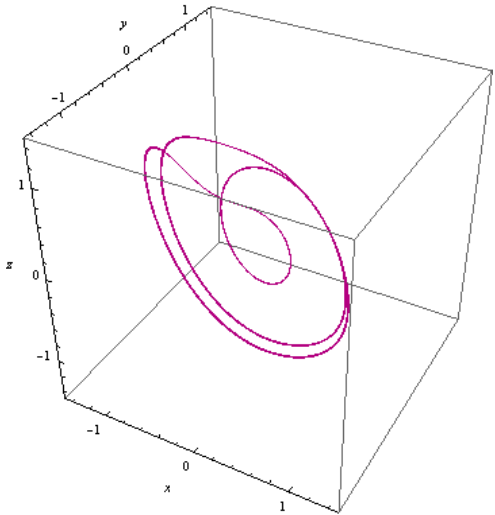 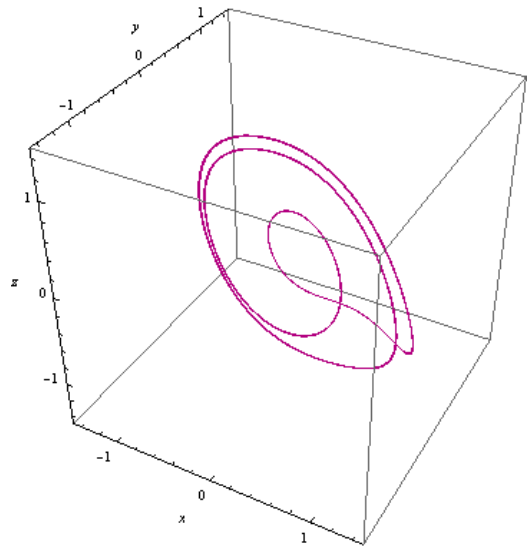

Figure 36. One limit closed orbit in the case $b = 0.3237$

Figure 37. Another limit closed orbit in the case $b = 0.3237$

There should be a series critical values, $b_{3,n}$, $n = 0, 1, 2, \cdots$, of the parameter $b$, distributed in the interval $(b_{3,inf}, b_3)$ satisfying

$$b_{3,inf} < \cdots < b_{3,n} < \cdots < b_{3,1} < b_{3,0} < b_3$$

such that

(3.1)  when $b_{3,0} \leq b < b_3$, the system (1) has a unique limit closed orbit of rotation number 3;

(3.2)  when $b_{3,n} \leq b < b_{3,n-1}$, the system (1) has a pair of spatial limit closed orbits of rotation number $2^n \times 3$;
(3.3)     $b_{3,inf} = \lim_{n \to \infty} b_{3,n}$.

By the similar numerical tests for the *series 1* and *series 2*,  the following two critical values are obtained as follows

$$b_{3,0} \approx 0.3184, \ b_{3,1} \approx 0.3177$$

But, it is very difficult to calculate the other critical values, limited by the numerical method.

However, based on the theory on Feigenbaum constant, we may obtain two estimations

$$b_{3,2} \approx b_{3,1} - \frac{b_{3,0} - b_{3,1}}{\delta} \approx 0.31755 \tag{8}$$

$$b_{3,inf} \approx b_{3,0} - \frac{b_{3,0} - b_{3,1}}{(1 - \frac{1}{\delta})} \approx 0.3175 \tag{9}$$

Figure 38 shows the rotation number of the limit closed orbit is still 12 in the case $b = 0.31755$

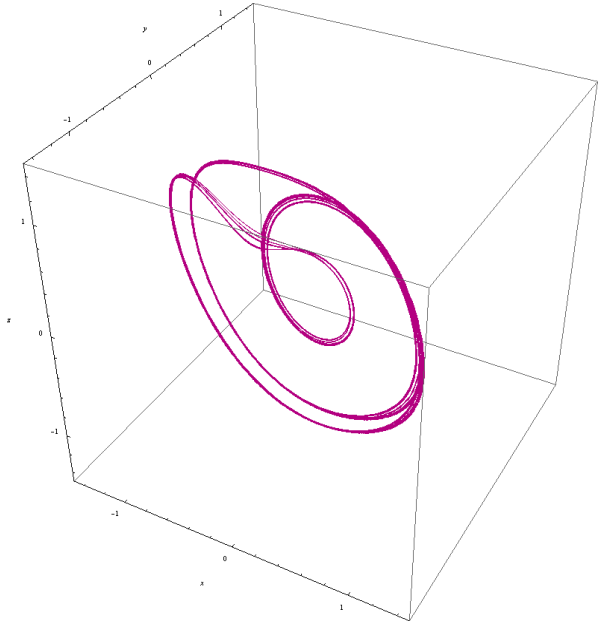

Figure 38. The limit closed orbit in the case $b = 0.31755$

Figure 39 shows the rotation number of the limit closed orbit has been doubled to 24 in the case $b = 0.31754$

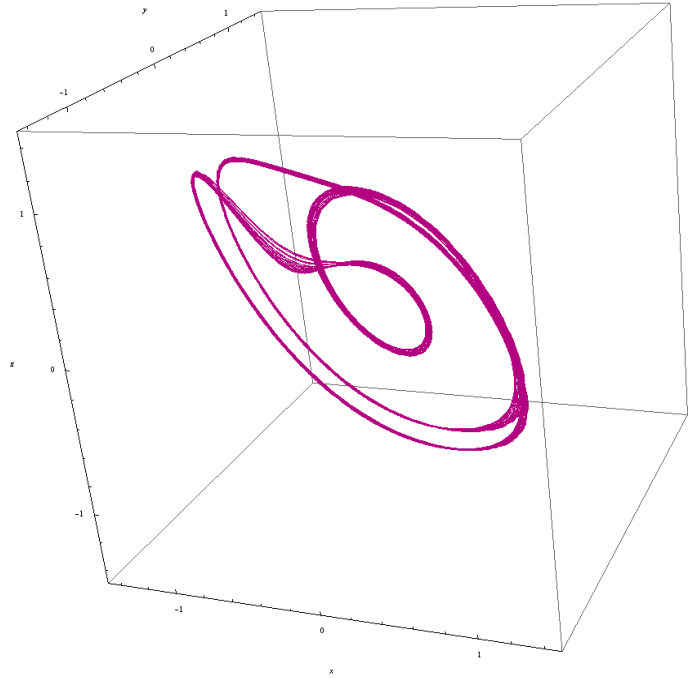

Figure 38. The limit closed orbit in the case $b = 0.31754$

The numerical test suggests that the set $A_{class}$ might still be formed with a pair of closed limit orbit in the case $b = 0.3176$ (see Figure 40), and that the set $A_{class}$ should not be formed with limit closed orbits in the case $b = 0.3175$ (see Figure 41).

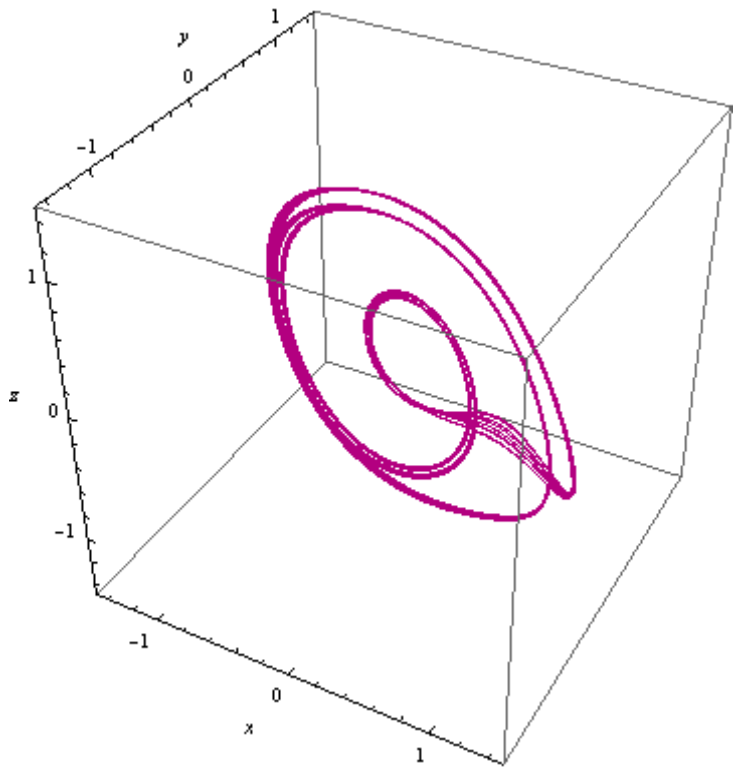

Figure 40. One branch of the set $A_{class}$ in the case $b = 0.31751$

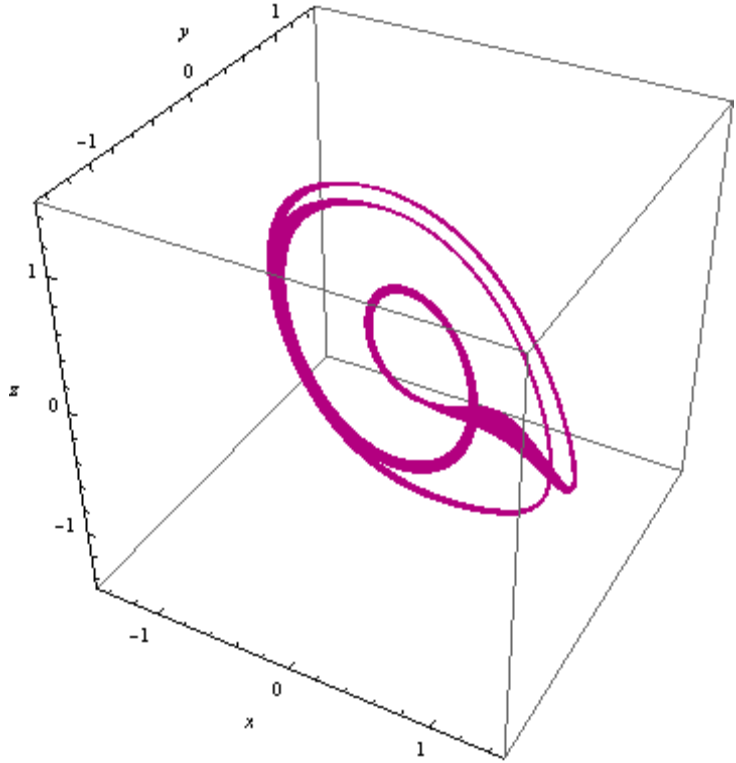

Figure 41.  One branch of the set  $A_{class}$  in the case  $b = 0.3175$

Again, we found that the estimations (8) and (9) are surprisingly effective. The universality of the Feigenbaum constant has been proven for the system (1).

## 4.  Some Further Discussion

Clearly, besides the discussed three series of the period-doubling cascades, there are still more complicated phenomena to be studied, for instance, in the interval $[b_3, b_{13,inf}]$ , there exist still other possible parameters, such that the system (1) has closed limit orbit.

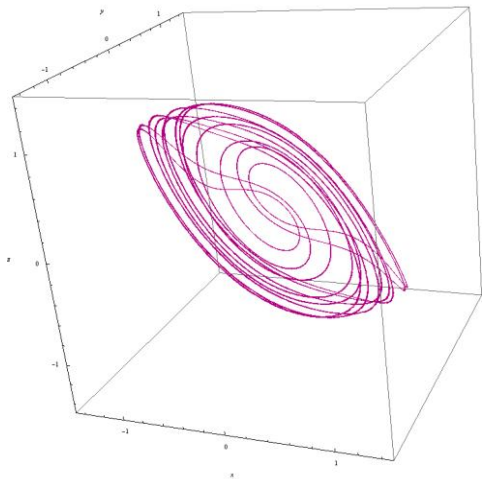

Figure 41.  The closed limit orbit in the case  $b = 0.3331$

For an instance, when b = 0.3331, the system (1) has a closed limit orbit with rotation number 13 (see Figure 42).

So it is in need to have more detailed research for the whole region of the parameter. It is quite possible to find more interesting phenomena.

It is a reasonable consideration that most of the bifurcation phenomena for the system (1) can be explained with the theory of the period-doubling cascades for the abstract dynamical systems given by Evelyn Sander and James A. Yorke. In fact, for a given closed limit orbit, one can define the Poincaré (or the first return) map of the system (1) (ref. [4]), it may play the role of the smooth one-parameter families of maps used in the abstract dynamical system theory.

However, this paper shows that there are still some particular bifurcation phenomena belonging to the ODE system (1), for instance, the dependence of the cascades on the topological character of rotation number of the limit closed orbit, This property has not been described in the present abstract dynamical system theory.

Another important issue should be mentioned here. Since the system (1) is symmetrical about the origin, sometimes，there are two symmetrical directions to go for the integral curves in the region $\mathfrak{M}$ between the two stable manifolds $M_{2d}^s(p_1)$ and $M_{2d}^s(p_2)$, for instances, the two symmetrical limit closed orbits, or the two pseudo "gaps" located respectively on the two manifolds (see ref. [1]). By the continuous dependence of the integral curve on its initial values, the region $\mathfrak{M}$ between the two manifolds should be separated into two open sub-regions and a non-empty closed sub set $S$, which is an invariant set of the system. Since this set plays the role to separate the region $\mathfrak{M}$, it will be called a separator of the system (1).

On the separator $S$, an integral curve cannot approach to any one of the above-mentioned two directions. Where would the integral curves on $S$ approach to? It is a reasonable conclusion that there should be a pseudo-attractor $A_{pseudo}$ on $S$, such that the integral curves on $S$ approach to it. In the case that the system (1) has a pair of symmetrical limit closed orbits, the pseudo-attractor may be formed with some closed cycle(s). And in the case that there exist two pseudo gaps, the pseudo-attractor $A_{pseudo}$ may be just the so-called "faint attractor", which was discussed in [1].

It is quite possible that the measure of the set $S$ is just zero, so that the chance to find numerically the pseudo-attractor $A_{pseudo}$ is zero. The more detailed discussion on the separator $S$, the pseudo-attractor $A_{pseudo}$ and faint attractor will be discussed in other papers.

The facts established for the system (1) have shown that this system is just an ideal model of the three-dimensional autonomous ODE system in the research of the bifurcation and chaos theory.